\documentclass[11pt]{amsart}

\usepackage{graphicx}
\usepackage{amscd}
\usepackage{amsmath}
\usepackage{amsfonts}
\usepackage{amssymb}
\usepackage{amsxtra}
\usepackage{xspace}
\usepackage{array}
\usepackage{cite}
\usepackage[all]{xy}
\xyoption{arc}

\theoremstyle{plain}
\newtheorem{theorem}{Theorem}
\newtheorem{lemma}{Lemma}
\newtheorem{corollary}{Corollary}
\newtheorem{proposition}{Proposition}

\theoremstyle{definition}
\newtheorem{definition}{Definition}

\numberwithin{equation}{section}

\newcommand{\be}{\begin{enumerate}}
\newcommand{\ee}{\end{enumerate}}
\newcommand{\beq}{\begin{equation}}
\newcommand{\eeq}{\end{equation}}
\newcommand{\bprop}{\begin{proposition}}
\newcommand{\eprop}{\end{proposition}}

\DeclareMathOperator{\ad}{ad}

\newcommand{\cleq}{\preccurlyeq}

\newcommand{\integers}{\mathbb{Z}}

\newcommand{\hkma}{hyperbolic Kac-Moody algebra\xspace}
\newcommand{\hkmas}{hyperbolic Kac-Moody algebras\xspace}
\newcommand{\gcm}{generalized Cartan matrix\xspace}
\newcommand{\kma}{\mathfrak{g}}
\newcommand{\s}[1][]{\ensuremath{\mathfrak{sl}_{#1}}}
\newcommand{\csa}{\mathfrak{h}}
\newcommand{\al}{\alpha}
\newcommand{\alno}{\alpha_{-1}}
\newcommand{\alc}{\Check{\alpha}}
\newcommand{\et}{\ensuremath{E_{10}}\xspace}
\newcommand{\rts}{\Phi}
\newcommand{\rtsp}{\Phi^+}
\newcommand{\rtsn}{\Phi^-}
\newcommand{\rtspre}{\rtsp_{re}}

\newcommand{\rtsre}{\rts_{re}}
\newcommand{\rtsim}{\rts_{im}}
\newcommand{\W}{W}

\newcommand{\rsd}{root subdiagram\xspace}
\newcommand{\rsds}{root subdiagrams\xspace}
\newcommand{\bp}{B$\mathbf{{}^\prime}$}

\begin{document}
\title{Embeddings of hyperbolic Kac-Moody algebras into $\mathbf{E_{10}}$}
\author{Sankaran Viswanath}
\address{Department of Mathematics\\
Penn State University\\
University Park, PA 16802, USA}
\email{viswanat@math.psu.edu}
\begin{abstract}
We show that the rank $10$ hyperbolic Kac-Moody algebra $E_{10}$ contains every  simply laced hyperbolic Kac-Moody algebra as a Lie subalgebra.
Our method is based on an extension of earlier work of Feingold and Nicolai. 
\end{abstract}

\maketitle

\section{Introduction}
Since their discovery \cite{kac-orig,moody-orig} by Kac and Moody, Kac-Moody algebras have been playing increasingly important roles in diverse subfields of mathematics and physics. The affine Kac-Moody algebras are by now as well understood as the finite dimensional simple Lie algebras classified by Cartan and Killing. 
Indefinite type Kac-Moody algebras however remain a notoriously intractable part of the theory. In spite of much work in this direction (see \cite{ff,kmw} and the references in \cite{fn}), obtaining detailed information about the structure of these Lie algebras seems out of reach at present. Most of the available results concern the subclass of {\em hyperbolic} Kac-Moody algebras. Such algebras only exist in ranks $2$-$10$ and can be completely classified \cite{sacl}  (see also \cite{extradigs} for some missing diagrams). Among these, the algebra \et has been singled out for its relevance to string theory and has received much attention in recent times (e.g: \cite{string1,string2, string3, string4, string5, string6, string7}).

In \cite{fn}, Feingold and Nicolai studied subalgebras of \hkmas. They
 showed that the rank $3$ \hkma $\mathcal{F}  = HA_1^{(1)}$ contains every rank $2$ \hkma with symmetric \gcm; in fact $\mathcal{F}$ was also shown to contain an infinite series of indefinite type Kac-Moody algebras. Analogously, it was shown that there are infinitely many inequivalent  Kac-Moody algebras of indefinite type that occur as Lie subalgebras of $E_{10}$.

This motivates the main question of this letter: {\em which Kac-Moody algebras of hyperbolic type occur as Lie subalgebras of \et ?} To answer this question, 
 we extend the method of  Feingold-Nicolai and formulate some general principles for constructing Lie subalgebras. These principles are then used to prove our main result: {\em Every} simply laced \hkma occurs as a Lie subalgebra of \et.
This statement can be viewed as further evidence of the distinguished role played by $E_{10}$ in the family of simply laced \hkmas. Though these subalgebras, being of hyperbolic type themselves, are  poorly understood objects, their embeddings in \et may neverthless shed more light on the structure of \et.

\medskip
\noindent
{\em Acknowledgements:} The author would like to thank Alex Feingold and Hermann Nicolai for their comments on an earlier draft of this letter. The author is also grateful to the anonymous referees for their careful reading of the manuscript and for their helpful suggestions that have led to a considerably improved exposition.
\section{Simply laced hyperbolic diagrams}
We start with a brief summary of the basic definitions and notation concerning
Kac-Moody algebras. The reader is referred to \cite{kac} for full
details. An $n \times n$ integer matrix $A = (a_{ij})$ is called a {\em generalized Cartan matrix} if it 
satisfies (i) $a_{ii} = 2 \; \forall i $, (ii) $a_{ij} \leq 0 \; \forall i \neq j$, and 
(iii)  $a_{ij} =0 \implies a_{ji} =0$. Given a generalized Cartan matrix $A$,
 one defines the associated Kac-Moody algebra $\kma = \kma(A)$ (of {\em rank} $n$)
 to be the Lie algebra with $3n$ generators $\alc_i, e_i, f_i \; (i=1, \cdots, n)$ subject to the relations $[\alc_i, \alc_j]=0, [\alc_i, e_j] = a_{ij} e_j, [\alc_i, f_j] = - a_{ij} f_j, [e_i, f_j] = \delta_{ij} \alc_i$ for all $i,j$, and the Serre relations $(\ad e_i)^{-a_{ij} + 1}(e_j) = 0 = (\ad f_i)^{-a_{ij} + 1}(f_j)$ for $i \neq j$.
If $A$ is nonsingular, the Cartan subalgebra $\csa$ is the span of the $\alc_i$; if $A$ is singular, one must also include certain derivations in $\csa$.
The simple roots $\al_j,\, j=1 \cdots n$ are linearly independent elements in $\csa^*$ which satisfy the relation $\al_j (\alc_i) = a_{ij} \; \forall i,j$; the integer lattice $\integers(\alpha_1, \cdots, \al_n)$ spanned by them is the {\em root lattice} $Q$. One has the root space decomposition:
$$ \kma = \csa \oplus (\oplus_{\alpha \in \rts} \; \kma_\alpha)$$
where each $\kma_\alpha =\{x \in \kma: [h,x]=\alpha(h) \, x, \, \forall h \in \csa\}$ is 
finite dimensional and $\rts = \{ \alpha \in \csa^*: \kma_\alpha \neq 0\} \subset Q$ is the set of roots of $\kma$. With $Q^+ = \integers_{\geq 0}(\alpha_1, \cdots, \alpha_n), \, \rtsp = \rts \cap Q^+$ and $\rtsn = - \rts^+$, we also have $\rts = \rtsp \cup \rtsn$. 
The Weyl group $\W \subset \mathrm{GL}(\csa^*)$ of $\kma$ is generated by the simple reflections $r_i$ of $\csa^*$ defined by $r_i (\al_j) =  \al_j - a_{ij} \al_i \; \forall j$.
The elements of $\rts$ that are $W$-conjugate to a simple root are termed {\em real}
 roots; the remaining are called {\em imaginary} roots. Letting $\rtsre$ and $\rtsim$ denote the subsets of real and imaginary roots respectively, we have $\rts = \rtsre \cup \rtsim$. 

If the matrix $A$ is  symmetric, $\kma$ is termed a {\em simply laced} Kac-Moody algebra. We will only consider simply laced algebras in the rest of this letter. In this case, 
the bilinear form $(,)$ on $\csa^*$ defined by $(\al_i, \al_j) = a_{ij} \;\forall i,j $  is symmetric and $\W$-invariant (more generally, such a form exists when $A$ is symmetrizable, i.e, when there exists a diagonal matrix $D$ such that $DA$ is symmetric). As mentioned above, the real roots are $\W$-conjugates of the simple roots; since $(\al_i, \al_i)=2 \; \forall i$, the $W$-invariance of the form means $(\beta, \beta)=2$ for all real roots $\beta$. For $\gamma \in \rtsim$, it is true that  $(\gamma, \gamma) \leq 0$ \cite{kac}. The {\em Dynkin diagram}
 associated to a symmetric generalized Cartan matrix $A$ is a graph on $n$ vertices, with vertices $i$ and $j$ joined by $|a_{ij}|$ edges for all $i \neq j$. If $|a_{ij}| > 1$, one also draws arrows on the edges joining vertices $i$ and $j$, pointing in the directions of both vertices (see for example Tables \ref{regs}, \ref{irregs}). The matrix $A$ (or equivalently its Dynkin diagram) is of {\em finite type} if all its principal minors are positive definite, of {\em affine type} if $\det A =0$ and all its proper principal minors are positive definite,  and of {\em indefinite type} otherwise. 

Let $X$ be an indefinite type Dynkin diagram. If every proper, connected subdiagram of $X$ is of finite or affine type, then $X$ (or alternately, its associated Kac-Moody algebra) is said to be 
of {\em hyperbolic} type. It is clear that $X$ is of hyperbolic type if and only if the deletion of any one vertex of $X$ results in a diagram each of whose connected components is of finite or affine type. 
Hyperbolic Kac-Moody algebras only exist in ranks $\leq 10$. 
For each $a \geq 3$, $\; \kma(\left[ \begin{smallmatrix} 2 & -a \\ -a &2 \end{smallmatrix}\right]) $ is a simply laced, 
hyperbolic Kac-Moody algebra of rank 2. In ranks $3$-$10$, there are exactly 23 simply laced hyperbolic Kac-Moody algebras  \cite{sacl}. Their Dynkin diagrams can be organized into 3 finite families (containing 15 diagrams, see Table \ref{regs}) and 8 exceptions (termed {\em irregular diagrams}, Table \ref{irregs}). We remark that if $X$ is a simply laced affine diagram, the notation $HX$ refers to its ``hyperbolic extension'', obtained by adding one more vertex, connected only to the zeroth node of $X$. We use $E_{10}$ as alternative notation for $H E_8^{(1)}$; we will also let $E_{10}$ denote 
the Kac-Moody algebra whose Dynkin diagram is $E_{10}$.
\begin{table}[t]
\begin{tabular}{l c}
$HA_1^{(1)}$ &
\def\objectstyle{\scriptscriptstyle}
$\xy
(0,0)*{\bullet} = "A" - (0,2)*{-1}; (6,0)*{\bullet} = "B"- (0,2)*{0};
(12,0)*{\bullet} = "C" - (0,2)*{1};
{\ar@2{<->}"B"; "C"};
 "A"; "B" **\dir{-};
\endxy
$\\\\
$HA_k^{(1)} \;\;(2 \leq k \leq 7)$ &
\def\objectstyle{\scriptscriptstyle}
$\xy 0;/r.2pc/:
(0,0)*{\bullet} = "A" - (0,2)*{1}; (6,0)*{\bullet} = "B"- (0,2)*{2};
 "A"; "B" **\dir{-};
(18,0)*{\bullet} = "C" - (0,2)*{k-1};
(24,0)*{\bullet} = "D" - (0,2)*{k};
"B"; "C" **\dir{.};
"D"; "C" **\dir{-};
(12,6)*{\bullet} = "E" + (2,1)*{0};
(12,10)*{\bullet} = "F" + (3,0)*{-1};
"A"; "E" **\dir{-};
"D"; "E" **\dir{-};
"E"; "F" **\dir{-};
\endxy$\\\\
$HD_k^{(1)} \;\;(4 \leq k \leq 8)$ &
\def\objectstyle{\scriptscriptstyle}
$\xy 0;/r.2pc/: 
(0,0)*{\bullet} = "A" - (0,2)*{-1}; (12,0)*{\bullet} = "B"- (0,2)*{2};
 "A"; "B" **\dir{-};
(12,4)*{\bullet} = "C" + (2,1)*{1};
"C"; "B" **\dir{-};
(6,0)*{\bullet} = "D"- (0,2)*{0};
(24,0)*{\bullet} = "E" - (0,2)*{}; (30,0)*{\bullet} = "F"- (0,2)*{k-1};
 "E"; "F" **\dir{-};
(24,4)*{\bullet} = "G" + (2,1)*{k};
"E"; "G" **\dir{-};
"E"; "B" **\dir{.};
\endxy$ \\\\
$HE_6^{(1)}$ &
\def\objectstyle{\scriptscriptstyle}
$\xy 0;/r.2pc/: 
(0,0)*{\bullet} = "A" - (0,2)*{-1}; 
(6,0)*{\bullet} - (0,2)*{0}; 
(12,0)*{\bullet} - (0,2)*{1}; 
(18,0)*{\bullet} = "B" - (0,2)*{2}; 
(24,0)*{\bullet} - (0,2)*{3}; 
(30,0)*{\bullet} = "C" - (0,2)*{4}; 
(18,4)*{\bullet} + (2,1)*{5};
(18,8)*{\bullet} = "D" + (2,1)*{6};
"A"; "C" **\dir{-};
"B"; "D" **\dir{-};
\endxy$ \\\\
$HE_7^{(1)}$ &
\def\objectstyle{\scriptscriptstyle}
$\xy 0;/r.2pc/: 
(-6,0)*{\bullet} = "A" - (0,2)*{-1}; 
(0,0)*{\bullet} - (0,2)*{0}; 
(6,0)*{\bullet} - (0,2)*{1}; 
(12,0)*{\bullet} - (0,2)*{2}; 
(18,0)*{\bullet} = "B" - (0,2)*{3}; 
(24,0)*{\bullet} - (0,2)*{4}; 
(30,0)*{\bullet} - (0,2)*{5}; 
(36,0)*{\bullet} = "C" - (0,2)*{6}; 
(18,4)*{\bullet} = "D" +  (2,1)*{7};
"A"; "C" **\dir{-};
"B"; "D" **\dir{-};
\endxy$ \\\\
$HE_8^{(1)} = E_{10}$ &
\def\objectstyle{\scriptscriptstyle}
$\xy 0;/r.2pc/: 
(-6,0)*{\bullet} = "A" - (0,2)*{-1}; 
(0,0)*{\bullet} - (0,2)*{0}; 
(6,0)*{\bullet} - (0,2)*{1}; 
(12,0)*{\bullet} - (0,2)*{2}; 
(18,0)*{\bullet} - (0,2)*{3}; 
(24,0)*{\bullet} - (0,2)*{4}; 
(30,0)*{\bullet} = "B" - (0,2)*{5}; 
(36,0)*{\bullet}  - (0,2)*{6}; 
(42,0)*{\bullet} = "C" - (0,2)*{7}; 
(30,4)*{\bullet} = "D" +  (2,1)*{8};
"A"; "C" **\dir{-};
"B"; "D" **\dir{-};
\endxy$ \\\\
 & 
\end{tabular}
\caption{``Regular'' simply laced hyperbolic diagrams}
\label{regs}
\end{table}

\begin{table}[t]
\begin{tabular}{l c}
$X_6$ &
\def\objectstyle{\scriptscriptstyle}
$\xy 
0;/r.2pc/:
(6,0)*{\bullet} = "A" + (0,2)*{2}; 
(0,6)*{\bullet} = "B" - (3,0)*{-1}; 
(0,-6)*{\bullet} = "C" - (2,0)*{0}; 
(12,6)*{\bullet} = "D" + (2,0)*{4}; 
(12,-6)*{\bullet} = "E" + (2,0)*{3}; 
(6,-6)*{\bullet} = "F" - (0,2)*{1}; 
"B"; "E" **\dir{-};
"C"; "D" **\dir{-};
"A"; "F" **\dir{-};
\endxy$\\\\
$Y_5$ &
\def\objectstyle{\scriptscriptstyle}
$\xy 
0;/r.2pc/:
(0,0)*{\bullet} = "D" - (0,2)*{1}; 
(6,0)*{\bullet} = "A" - (0,2)*{0}; 
(0,6)*{\bullet} = "B" - (2,0)*{2}; 
(6,6)*{\bullet} = "C" + (2,0)*{3}; 
(12,0)*{\bullet} = "E" - (0,2)*{-1}; 
"D"; "B" **\dir{-};
"B"; "C" **\dir{-};
"C"; "A" **\dir{-};
"D"; "E" **\dir{-};
"E"; "B" **\crv{(10,15)};
\endxy$\\\\

$Y_4$ &
\def\objectstyle{\scriptscriptstyle}
$\xy 
0;/r.2pc/:
(0,0)*{\bullet} = "D" - (2,0)*{1}; 
(0,6)*{\bullet} = "B" - (2,0)*{2}; 
(6,3)*{\bullet} = "C" - (0,2)*{0}; 
(12,3)*{\bullet} = "E" - (0,2)*{-1}; 
"D"; "B" **\dir{-};
"B"; "C" **\dir{-};
"C"; "D" **\dir{-};
"C"; "E" **\dir{-};
"E"; "B" **\crv{(8,8)};
\endxy$\\\\

$Z_4$ &
\def\objectstyle{\scriptscriptstyle}
$\xy 
0;/r.2pc/:
(0,0)*{\bullet} = "D" - (2,0)*{1}; 
(0,6)*{\bullet} = "B" - (2,0)*{2}; 
(6,3)*{\bullet} = "C" - (0,2)*{0}; 
(12,3)*{\bullet} = "E" + (3,0)*{-1}; 
"D"; "B" **\dir{-};
"B"; "C" **\dir{-};
"C"; "D" **\dir{-};
"C"; "E" **\dir{-};
"E"; "B" **\crv{(8,8)};
"E"; "D" **\crv{(8,-2)};
\endxy$\\\\

$Y_3$ &
\def\objectstyle{\scriptscriptstyle}
$\xy
(-1,0)*{\bullet} = "A" - (2,0)*{1}; 
(7,0)*{\bullet} = "B" + (2,0)*{0};
(3,-6)*{\bullet} = "C" - (0,2)*{-1};
{\ar@2{<->}"B"; "A"};
 "C"; "B" **\dir{-};
 "C"; "A" **\dir{-};
\endxy
$\\\\

$T_2$ &
\def\objectstyle{\scriptscriptstyle}
$\xy
(-1,0)*{\bullet} = "A" - (2,0)*{1}; 
(7,0)*{\bullet} = "B" + (2,0)*{0};
(3,-6)*{\bullet} = "C" - (0,2)*{-1};
{\ar@2{<->}"B"; "A"};
{\ar@2{<->}"C"; "A"};
{\ar@2{<->}"B"; "C"};
\endxy
$\\\\
$T_1$ &
\def\objectstyle{\scriptscriptstyle}
$\xy
(-1,0)*{\bullet} = "A" - (2,0)*{1}; 
(7,0)*{\bullet} = "B" + (2,0)*{0};
(3,-6)*{\bullet} = "C" - (0,2)*{-1};
{\ar@2{<->}"B"; "A"};
{\ar@2{<->}"B"; "C"};
"C"; "A" **\dir{-};
\endxy
$\\\\
$T_0$ &
\def\objectstyle{\scriptscriptstyle}
$\xy
(-1,0)*{\bullet} = "A" - (2,0)*{1}; 
(7,0)*{\bullet} = "B" + (2,0)*{0};
(3,-6)*{\bullet} = "C" - (0,2)*{-1};
{\ar@2{<->}"B"; "A"};
{\ar@2{<->}"B"; "C"};
\endxy
$\\\\
\end{tabular}
\caption{Irregular simply laced hyperbolics}
\label{irregs}
\end{table}

The main result of this letter is the following :
\begin{theorem}\label{mainprop}
Given any simply laced hyperbolic Kac-Moody algebra $\kma$, there is a Lie subalgebra of $E_{10}$ that is isomorphic to $\kma$.
\end{theorem}

To prove this theorem, we will use the following result \cite[Theorem 3.1]{fn} (restated for our situation):

\begin{theorem}\label{fnthm}
(Feingold-Nicolai) Let $A$ be an $n \times n$  symmetric generalized Cartan matrix and $\kma = \kma(A)$. Suppose $k \leq n$ and 
$\beta_1, \cdots, \beta_k \in \rtspre$ such that $\beta_i - \beta_j \not\in \rts \; \forall i \neq j$. Choose nonzero elements $E_{\beta_i} \in \kma_{\beta_i}$ and $F_{\beta_i} \in \kma_{-\beta_i}$ for all $i$. The Lie subalgebra of $\kma$ generated by $\{ E_{\beta_i}, F_{\beta_i}: i=1, \cdots, k\}$ is a Kac-Moody algebra of rank $k$, with generalized Cartan matrix $C = [c_{ij}] = [(\beta_i, \beta_j)]$.
\end{theorem}

We also state the following lemma (implicit in \cite{fn}) that simplifies the task of verifying the hypothesis of theorem \ref{fnthm}.
\begin{lemma}
With notation as in theorem \ref{fnthm}, suppose $\beta, \gamma \in \rtspre$ satisfy $(\beta, \gamma) \leq 0$, then $\beta - \gamma \not\in \rts$. 
\end{lemma}

\noindent
{\em Proof:} Since $(\beta, \beta) = (\gamma, \gamma) =2$, we have $ (\beta - \gamma, \beta - \gamma) = (\beta, \beta) +  (\gamma, \gamma) - 2(\beta, \gamma) \geq 4$. Since $\beta - \gamma$ has 
positive norm, it cannot be in $\rtsim$; since its norm is not 2, it cannot be in $\rtsre$ either. 

Thus, in theorem \ref{fnthm}, if $\beta_1, \cdots, \beta_k \in \rtspre$ are such that $ C:=[(\beta_i, \beta_j)]$ is a generalized Cartan matrix (i.e, $(\beta_i, \beta_j) \leq 0 \; \forall i \neq j$), the subalgebra generated by the $E_{\beta_i}, F_{\beta_i} $ is automatically a Kac-Moody algebra with \gcm $C$.

\begin{definition} \label{rsddef}
In the notation of theorem \ref{fnthm}, if $D$ denotes the Dynkin diagram of $\kma$ and $D^\prime$ is the Dynkin diagram corresponding to the generalized Cartan matrix $C = [(\beta_i, \beta_j)]$, we will call $D^\prime$ a 
{\bf{root subdiagram}} of $D$. We will denote this by $D^\prime \cleq D$.
\end{definition}

As is clear from the definition, root subdiagrams need not be subdiagrams of $D$.
The following statement is easily seen to imply theorem \ref{mainprop}.
\begin{proposition} \label{mainpropprime}
Every simply laced hyperbolic Dynkin diagram occurs as a \rsd of $E_{10}$.
\end{proposition}

Sections \ref{three} and \ref{four} will be devoted to the proof of proposition \ref{mainpropprime}.

\section{Constructing root subdiagrams}\label{three}
The goal of this section is to formulate four general principles for constructing root subdiagrams of certain kinds of Dynkin diagrams. For the first two principles, we let $X$ denote a simply laced {\em affine}
 Dynkin diagram and $Y=HX$, the hyperbolic extension of $X$ :
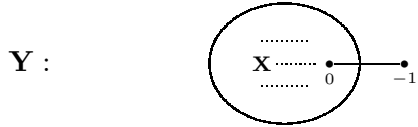
$$
{\mathbf Y}: \hspace{2cm}
\def\objectstyle{\scriptscriptstyle}
\def\labelstyle{\scriptscriptstyle}
\xy
(0,0)*{\xycircle(10,8){-}};
(6,0)*{\bullet} = "A" - (0,2)*{{0}}; 
(16,0)*{\bullet} = "B" - (0,2)*{{-1}};
"A"; "B" **\dir{-};
(-3,0)*{{\bf \scriptstyle X}};
(4,0)*{}; (-1,0)*{} **\dir{.};
(3,3)*{}; (-3,3)*{} **\dir{.};
(3,-3)*{}; (-3,-3)*{} **\dir{.};
\endxy
$$

Let the simple roots of $X$ be $\al_0, \al_1, \cdots, \al_r$ and let $\alno$ be the simple root corresponding to the extra vertex of $Y$. Let $\delta_X$ denote the null root of $X$; thus  
$\delta_X = \sum_{i=0}^r n_i \al_i$ ($n_i \in \integers_{\geq 0}, \, n_0=1$) and 
$(\delta_X ,\al_i)=0 \; \forall i=0, \cdots, r$. Note that this means $(\delta_X,\alpha_{-1})=-1$.

\medskip
\noindent
{\bf Principle A:} Define $\beta_{-1}, \beta_0, \cdots, \beta_r$ as follows:
$$ \beta_i = \al_i \;\;\;\; (0 \leq i \leq r),  \;\;\;\;\;\; \beta_{-1} = r_{-1} (\delta_X + \al_0)$$
where $r_{-1} \in \W$ is the simple reflection corresponding to $\alno$.  Since $\delta_X + \al_0 = 2\al_0 + 
\sum_{i=1}^r n_i \al_i$ and the vertex $-1$ is only connected to vertex $0$, one observes that 
$$ \beta_{-1} = \delta_X + \al_0 + 2 \alno$$
Since $\delta_X + \al_0$ is a real root of the affine Kac-Moody algebra $\kma(X)$, it is also a real root
of  $\kma(Y)$; thus its $\W$-conjugate, viz $\beta_{-1}$, is in $\rtspre$.

Next, the bilinear form $(,)$ : 
$$ (\beta_i, \beta_j) = (\al_i, \al_j) \;\;  (0 \leq i,j \leq r) $$
For $0 \leq j \leq r$, we also have :
\begin{align*}
(\beta_{-1}, \beta_j) &= ( \delta_X + \al_0 + 2 \alno, \al_j) = (\al_0 + 2 \alno, \al_j) \\
 &=\begin{cases} 0 & j \neq 0 \text{ and } j \text{ is not connected to vertex } 0 \\
 (\al_0,\al_j) & j \neq 0 \text{ and } j \text{ is  connected to vertex } 0 \\
0 & j=0
\end{cases}
\end{align*}
The $\beta_i$ thus satisfy the hypothesis of theorem \ref{fnthm}. 
Next consider the \rsd of $Y$ formed by the $\beta_{-1},\beta_0, \cdots, \beta_r$:  we observe that the diagram formed by its subset $\{\beta_i\}_{i=0}^r$ is just $X$; we then add an extra vertex corresponding to $\beta_{-1}$. This vertex is not connected to vertex $0$ anymore, but instead to all the neighbors of vertex $0$ in X.
For instance, applying principle {\bf A} to $Y = E_{10}, X=E_8^{(1)}$, the $\beta_i$ generate 
the following \rsd (upto renumbering of vertices, this is $HD_8^{(1)}$ from Table \ref{regs}) :
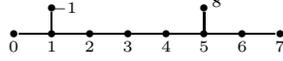
$$
\def\objectstyle{\scriptscriptstyle}
\xy 0;/r.2pc/: 
(6,4)*{\bullet} = "Y" + (2,0)*{-1}; 
(0,0)*{\bullet} = "A" - (0,2)*{0}; 
(6,0)*{\bullet} ="Z"  - (0,2)*{1}; 
(12,0)*{\bullet} - (0,2)*{2}; 
(18,0)*{\bullet} - (0,2)*{3}; 
(24,0)*{\bullet} - (0,2)*{4}; 
(30,0)*{\bullet} = "B" - (0,2)*{5}; 
(36,0)*{\bullet}  - (0,2)*{6}; 
(42,0)*{\bullet} = "C" - (0,2)*{7}; 
(30,4)*{\bullet} = "D" +  (2,1)*{8};
"A"; "C" **\dir{-};
"B"; "D" **\dir{-};
"Y"; "Z" **\dir{-};
\endxy$$

\begin{figure}
\end{figure}
Thus
\begin{equation} \label{hd8}
HD_8^{(1)} \cleq E_{10}
\end{equation}

\noindent
{\bf Principle B:} Let $X,Y, \, \al_i \;\; (-1 \leq i \leq r), \, \delta_X$ be as above. Choose
a vertex $1 \leq p \leq r$ of $X$, i.e,  a vertex of the finite type diagram underlying $X$.
Define
$$ \beta_i = \al_i \;\;\; (-1 \leq i \leq r, \, i \neq p), \;\;\;\;\; \beta_p = \al_p + \delta_X $$
Again, $\beta_i \in \rtspre \;\, \forall i=-1, \cdots, r$ and the bilinear form satisfies 
$$(\beta_i, \beta_j) = (\al_i, \al_j), \; \;  -1 \leq i,j \leq r; \;\; i,j \neq p$$
\begin{align*}
(\beta_p, \beta_i) &= (\al_p + \delta_X, \beta_i), \;\;\; -1 \leq i \leq r; \;i \neq p \\
&=\begin{cases} (\al_p, \al_i) & 0 \leq i \leq r; \; i \neq p \\
-1 & i=-1
\end{cases}
\end{align*}
The \rsd of $Y$ formed by the $\beta_{-1},\beta_0, \cdots, \beta_r$  
is thus the same as the original diagram $Y$ with one additional edge
connecting the vertex $-1$ to the chosen node $p$ of $X$.
For instance, figures \ref{ha8} and \ref{pt}  show the root subdiagrams $HA_8^{(1)}$ and $P_{10}$ 
obtained by  
applying Principle {\bf B} to $Y=E_{10}$ and $p=7$, $p=8$.
\begin{figure}
\def\objectstyle{\scriptscriptstyle}
$\xy
(-6,0)*{\bullet} = "A" - (0,3)*{-1}; 
(0,0)*{\bullet} - (0,1)+*{0}; 
(6,0)*{\bullet} - (0,2)*{1}; 
(12,0)*{\bullet} - (0,2)*{2}; 
(18,0)*{\bullet} - (0,2)*{3}; 
(24,0)*{\bullet} - (0,2)*{4}; 
(30,0)*{\bullet} = "B" - (0,2)*{5}; 
(36,0)*{\bullet}  - (0,2)*{6}; 
(42,0)*{\bullet} = "C" - (0,2)*{7}; 
(30,4)*{\bullet} = "D" +  (2,1)*{8};
"A"; "C" **\dir{-};
"B"; "D" **\dir{-};
"A"; "C" **\crv{(21,-16)};
\endxy$
\caption{$HA_8^{(1)}$}
\label{ha8}
\end{figure}

\begin{figure}
\def\objectstyle{\scriptscriptstyle}
$\xy
(-6,0)*{\bullet} = "A" - (0,2)*{-1}; 
(0,0)*{\bullet} - (0,2)*{0}; 
(6,0)*{\bullet} - (0,2)*{1}; 
(12,0)*{\bullet} - (0,2)*{2}; 
(18,0)*{\bullet} - (0,2)*{3}; 
(24,0)*{\bullet} - (0,2)*{4}; 
(30,0)*{\bullet} = "B" - (0,2)*{5}; 
(36,0)*{\bullet}  - (0,2)*{6}; 
(42,0)*{\bullet} = "C" - (0,2)*{7}; 
(30,4)*{\bullet} = "D" +  (2,1)*{8};
"A"; "C" **\dir{-};
"B"; "D" **\dir{-};
"A"; "D" **\crv{(15,8)};
\endxy
$
\caption{$P_{10}$}
\label{pt}
\end{figure}
Thus $HA_8^{(1)} \cleq E_{10}$ and $P_{10} \cleq E_{10}$; we remark that $HA_8^{(1)}$ and $P_{10}$ are indefinite type diagrams that are not hyperbolic.

\medskip
\noindent
{\bf Principle \bp:} In fact, one can pick any subset $F \subset \{1,2, \cdots, r\}$ and let 
$$ \beta_i = \al_i \; (-1 \leq i \leq r, i \not\in F), \;\;\;\;\; \beta_i = \al_i + \delta_X \; (i \in F) $$
A similar argument shows that the resulting root subdiagram is obtained by connecting vertex $-1$ to all 
the vertices in $F$.

\medskip
\noindent
{\bf  Principle C:} (shrinking) Suppose a Dynkin diagram $Z$ with $n$ vertices labeled $1 \cdots n$ 
has a subset  $\{k, k+1, \cdots, k+p-1\}$ of $p$ vertices which forms a subdiagram isomorphic to the finite type diagram $A_p$:
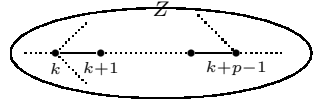
$$
\def\objectstyle{\scriptscriptstyle}
\xy
(14,0)*{\xycircle(20,6){-}};
(14,6)*{\scriptstyle Z};
(0,0)*{\bullet} = "A" - (0,2)*{k}; 
(6,0)*{\bullet} = "B"- (0,2)*{k+1};
(18,0)*{\bullet} = "C" ;
(24,0)*{\bullet} = "D" - (0,2)*{k+p-1};
 "A"; "B" **\dir{-};
"B"; "C" **\dir{.};
"D"; "C" **\dir{-};
"A"; "A" - (-4,4) **\dir{.};
"A"; "A" - (-4,-4) **\dir{.};
"A"; "A" - (4,0) **\dir{.};
"D"; "D" + (-5,5) **\dir{.};
"D"; "D" + (6,0) **\dir{.};
\endxy
$$
Letting $\al_i, \, i=1, \cdots,n $ denote the simple roots of $Z$, define
$$ \beta_j = \al_j, \;\; \text{ if } j<k \; \text{ or } j \geq k+p\, ; \;\;\;\;\;\;\;\;\; \bar{\beta} = \sum_{i = k}^{k+p-1} \al_i$$
The root subdiagram formed by $\{\beta_j: j<k\} \cup \{\bar{\beta} \} \cup \{\beta_j:  j \geq k+p \}$ has $p-1$ fewer vertices than $Z$; the vertex correponding to 
$\bar{\beta}$ is now connected to the original neighbors of vertex $k$ as well to those of vertex $k+p-1$, while the rest of the diagram remains unchanged. We note that if vertices $k$ and $k+p-1$ have a common neighbor $s$, then the vertex $\bar{\beta}$ is connected to $s$ by 2 lines (with arrows pointing both toward $s$ and $\bar{\beta}$) . An application of  Principle {\bf C} shows for instance that all members of the family $HA_k^{(1)} \; (1 \leq k \leq 7)$   occur as root subdiagrams of $HA_8^{(1)}$.
Similarly $HD_k^{(1)} \; (4\leq k \leq 8)$ occur as root subdiagrams of $HD_8^{(1)}$.

\medskip
\noindent
{\bf  Principle D:} (deletion) Given a Dynkin diagram $Z$ with $n$ vertices, deleting any subset of vertices 
 (and all incident edges) clearly gives us a root subdiagram of $Z$.

\section{Proof of main result}\label{four}
\noindent
The goal of this section is  to prove proposition \ref{mainpropprime} using these four principles.

\medskip
\noindent
{\em Proof of proposition \ref{mainpropprime}:} In \cite{fn}, it was shown that all the rank 2 simply laced hyperbolic Dynkin diagrams occur as \rsds of $HA_1^{(1)} = \mathcal{F}$. So, it is enough to show that the 23 diagrams in ranks 3-10 (Tables \ref{regs}, \ref{irregs}) occur as \rsds of $E_{10}$, since $HA_1^{(1)}$ occurs among these 23. Now, we had already observed that (i) $HD_8^{(1)} \cleq E_{10}$ (equation \eqref{hd8})
(ii) $HA_8^{(1)} \cleq E_{10}$ (figure \ref{ha8}) (iii) $P_{10} \cleq E_{10}$ (figure \ref{pt}).
 As remarked earlier, principle {\bf C} then implies that  $HA_k^{(1)}\;(1 \leq k \leq 7)$ and  $HD_k^{(1)} 
\; (4\leq k \leq 8)$ occur as root subdiagrams of $E_{10}$.

Next, we consider $HE_k^{(1)}, k=6,7,8$. When $k=8$, this is just $E_{10}$ itself. 
Consider $k=7$; 
we observe that a diagram isomorphic to $HE_7^{(1)}$ is obtained on deletion of 
 vertex $0$ from $HA_8^{(1)}$ 
(figure \ref{ha8}). Thus $HE_7^{(1)} \cleq HA_8^{(1)} \cleq E_{10}$.
Similarly, observe that $HE_6^{(1)} \cleq P_{10}$, and is 
obtained by deletion of the two vertices of $P_{10}$ numbered $0$ and $1$ in figure \ref{pt}.
So, the hyperbolic diagrams in the three families $HA_k^{(1)}, HD_k^{(1)}, HE_k^{(1)}$ are all \rsds of $E_{10}$.
For the 8 irregular diagrams (Table \ref{irregs}):
\be
\item $X_6$ : We start with $HD_4^{(1)} \cleq E_{10}$. Applying principle {\bf A} to $HD_4^{(1)}$ clearly 
gives us $X_6$. Thus $X_6 \cleq HD_4^{(1)} \cleq E_{10}$.
\item Diagrams $Y_k \; (k=3,4,5)$ are clearly obtained from $HA_{k-2}^{(1)}$ by 
application of principle {\bf B} (for suitable choice of $p$). For instance we have
$$
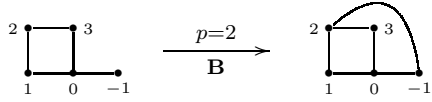
\def\objectstyle{\scriptscriptstyle}
\xy
(0,0)*{\bullet} = "D" - (0,2)*{1}; 
(6,0)*{\bullet} = "A" - (0,2)*{0}; 
(0,6)*{\bullet} = "B" - (2,0)*{2}; 
(6,6)*{\bullet} = "C" + (2,0)*{3}; 
(12,0)*{\bullet} = "E" - (0,2)*{-1}; 
"D"; "B" **\dir{-};
"B"; "C" **\dir{-};
"C"; "A" **\dir{-};
"D"; "E" **\dir{-};
(40,0)*{\bullet} = "D" - (0,2)*{1}; 
(46,0)*{\bullet} = "A" - (0,2)*{0}; 
(40,6)*{\bullet} = "B" - (2,0)*{2}; 
(46,6)*{\bullet} = "C" + (2,0)*{3}; 
(52,0)*{\bullet} = "E" - (0,2)*{-1}; 
"D"; "B" **\dir{-};
"B"; "C" **\dir{-};
"C"; "A" **\dir{-};
"D"; "E" **\dir{-};
"E"; "B" **\crv{(50,15)};
{\ar^{p=2}_{\bf B} (18,3)*{};(32,3)*{}};
\endxy$$

\item Diagram $Z_4$ is a \rsd of $HA_2^{(1)}$ via principle {\bf \bp}.
\item The remaining rank 3 diagrams $T_m \;(m=0,1,2)$ can each be obtained as a \rsd of $HA_1^{(1)}$ by choosing $\beta_{-1}, \, \beta_0, \, \beta_1$ as follows:
\begin{align*}
\beta_1 &= m\delta + \al_1 = (m+1) \al_1 + m \al_0\\
\beta_0 &= \delta + \al_0 =  \al_1 + 2 \al_0\\
\beta_{-1} &= \al_{-1}
\end{align*}
Here, the $\al_i$ are the simple roots of $HA_1^{(1)}$ and $\delta = \al_0 + \al_1$ is the null root of the underlying affine diagram $A_1^{(1)}$. We note that this is just a minor modification of theorem $3.3$
 of \cite{fn}. \qed
\ee

We also note the following fact: (notation as in definition \ref{rsddef}) if $D^\prime \cleq D$, then 
 the reflections $r_{\beta_i}$ generate  a subgroup of the Weyl group of  $D$  isomorphic to 
the Weyl group of $D^\prime$ \cite{dyer}.
Thus, as a corollary to proposition \ref{mainpropprime}, we infer that all simply laced Weyl groups of hyperbolic type occur as reflection subgroups (i.e subgroups generated by reflections) of the Weyl group of $E_{10}$.

\section{Disconnected root subdiagrams}
Thus far, we  were primarily concerned with the task of realizing {\em connected} simply laced hyperbolic Dynkin diagrams as \rsds of $E_{10}$. A much harder task is to identify {\em all} (possibly disconnected) simply laced \rsds of $E_{10}$. We conclude this letter with the following proposition which realizes two interesting disconnected Dynkin diagrams of rank $10$ as \rsds of $E_{10}$.

{\em Notation}: If $X,Y$ are two connected Dynkin diagrams, we will let $X \oplus Y$ denote their disjoint union (i.e, the disconnected diagram with components isomorphic to $X$ and $Y$).

\begin{proposition}\label{thmtwo}
The rank 10 diagrams  $HE_7^{(1)} \oplus A_1$ and $HE_6^{(1)} \oplus A_2$ occur as \rsds of $E_{10}$.
\end{proposition}
To prove this proposition, we will use the following characterization of real roots \cite[Proposition 5.10]{kac}.
\begin{theorem} \label{hyproot}
Let $\kma$ be a simply laced Kac-Moody algebra of finite, affine or hyperbolic type and let $Q$ be its root lattice. Then $$\rtsre = \{ \beta \in Q: (\beta,\beta)=2\}$$
\end{theorem}

\noindent
{\em Proof of proposition \ref{thmtwo}:} 
Let $\al_i \; (-1 \leq i \leq 8)$ be the simple roots of \et and $\delta$ be the null root of $E_8^{(1)}$. 
We recall that  $HE_7^{(1)} \cleq HA_8^{(1)} \cleq E_{10}$. Unravelling this chain of inclusions, we obtain the simple roots $\{\beta_i\}_{i=1}^9$ of  $HE_7^{(1)}$ to be 

$$\beta_i = \al_i \;(1 \leq i \leq 6), \;\;\; \beta_7 = \delta + \al_7, \;\;\;  \beta_8 = \al_8, \;\;\; 
\beta_9 = \al_{-1}$$
Note that we are now numbering the vertices of $HE_7^{(1)}$ with the labels $1\cdots9$ rather than $-1\cdots 7$.
If we can produce $\gamma \in \rtspre$ such that 
\begin{equation} \label{star}
(\gamma, \beta_i) =0  \text{ for all } 1 \leq i \leq 9 
\end{equation}
then the elements $\{\beta_i\}_{i=1}^9 \cup \{\gamma\}$ will generate $HE_7^{(1)} \oplus A_1$.
Let $\{\Lambda_i\}_{i=-1}^8$ denote the fundamental weights of \et i.e $(\Lambda_i, \al_j)=\delta_{ij} \; \forall i,j$ (see \cite{kmw} for a table of the $\Lambda_i$ as linear combinations of $\al_j$).
Equation \eqref{star} can be easily seen to imply that $\gamma \in \text{ span }(\Lambda_7 - 3 \Lambda_0)$. Another easy computation gives $(\Lambda_7 - 3 \Lambda_0, \Lambda_7 - 3 \Lambda_0) = 2$ and $\Lambda_7 - 3 \Lambda_0 \in \integers_{\geq 0} (\al_{-1}, \cdots, \al_8)$. Thus, using theorem \ref{hyproot}, $\gamma =\Lambda_7 - 3 \Lambda_0 $ is the required root.

Next, we have $HE_6^{(1)} \cleq P_{10} \cleq E_{10}$. The simple roots $\{\beta_i\}_{i=1}^8$  of 
$HE_6^{(1)}$ are now:
$$\beta_1 = \al_{-1}, \;\;\; \beta_i = \al_i \;\; (2 \leq i \leq 7), \;\;\; \beta_8 = \al_8 + \delta $$
As above, we look for all $\gamma \in \rtspre$ such that 
$(\gamma, \beta_i) =0 \; \forall i=1, \cdots, 8$. This condition implies 
$\gamma \in \text{ span }(\Lambda_8 - 2 \Lambda_1, \Lambda_1 - 2 \Lambda_0)$. Letting
$\gamma_1 = \Lambda_8 - 2 \Lambda_1$ and $\gamma_2 = \Lambda_1 - 2 \Lambda_0$, one computes
$(\gamma_j, \gamma_j)=2$ and $\gamma_j \in \integers_{\geq 0} (\al_{-1}, \cdots, \al_8)$ for $j=1,2$. 
Further $(\gamma_1, \gamma_2) = -1$. Thus $\{\beta_i\}_{i=1}^8 \cup \{\gamma_1, \gamma_2\}$ generates $HE_6^{(1)}
 \oplus A_2$. This completes the proof. \qed

\begin{corollary}
$E_{10}$ has Lie subalgebras isomorphic to $\kma(HE_6^{(1)}) \oplus \s[3] $ and  $\kma(HE_7^{(1)}) \oplus 
\s[2]$.
\end{corollary}

\end{document}